\documentclass[10pt,a4paper,reqno]{amsart}
\usepackage{amsfonts,amsthm,latexsym,amsmath,amssymb,amscd,amsmath, epsf}

\usepackage[T2A]{fontenc}
\usepackage[cp1251]{inputenc}
\usepackage{latexsym,amssymb,amsthm}
\usepackage{url,graphics} 
\usepackage[dvips]{graphicx,epsfig} 

\newtheorem{thm}{Theorem}

\def\q#1.{{\bf #1.}}
\renewcommand\geq{\geqslant}
\renewcommand\leq{\leqslant}

\newcommand{\be}{\begin{equation}}
\newcommand{\ee}{\end{equation}}

\begin{document}
          \numberwithin{equation}{section}

          \title[A  short proof of Kotzig's theorem]
          {A short proof of Kotzig's theorem}

\author[G. Nenashev]
{Gleb Nenashev}
\thanks{ \scriptsize This~research~is~supported~by~RF~Government~grant~11.G34.31.0026~and~by~JSC~"Gazprom~Neft"}
\address{  Chebyshev Laboratory, St. Petersburg State University, 14th Line, 29b, Saint Petersburg, 199178 Russia.}
\email{glebnen@mail.ru}

\begin{abstract}  
A new shortest proof of Kotzig's Theorem about  graphs with  unique  perfect matching is presented in this paper.
It is well known that Kotzig's theorem is a consequence of Yeo's Theorem about edge-colored graph without alternating cycle.
We present a proof of Yeo's Theorem based on the  same ideas as our proof of Kotzig's theorem.
 \end{abstract}

\maketitle

\section{\bf Introduction}

The well-known theorem of A.\,Kotzig was proved for the first time  in~\cite{K}.

\begin{thm}{\bf(A.\,Kotzig, 1959)} Let $G$ be a connected graph with unique perfect matching.  Then $G$ has a bridge that belongs to this matching. 
\label{Kot}
\end{thm}

However, the proof in~\cite{K} was tedious.  The shortest proof of Kotzig's theorem that is known now is to derive it from the following  theorem of A.\,Yeo~\cite{Y}.

We denote by $G-e$ graph $G$ without the edge $e$ and by $G-v$ graph $G$ without the vertex $v$ and all  edges incident to it.  

\begin{thm}{\bf(A.\,Yeo, 1997)} 
All edges of a graph $G$ are colored such that there is no alternating cycles (i.e. each cycle has two adjacent edges of the same colors). 
Then  $G$ contains a vertex $v$ such that every connected component of $G-v$ is joined to $v$ with edges of  one color.
\label{Ye}
\end{thm}

This theorem have rather short and elegant proof using the method of alternating chains. Let us also mention, that a particular case of Yeo's Theorem for coloring with two colors was proved by Grossman and Haggkvist in 1982~\cite{GH}.

Our short  proof of the Theorem~\ref{Kot} is based on analyzing of a minimal counterexample. 
Also we show that our method works in  Theorem~\ref{Ye}.

\medskip
\noindent
{\bf Acknowledgement.} The author is grateful to Dmitri Karpov for his comments and for discussion of proofs.

\section{\bf Proofs}

\begin{proof}[\bf Proof of Theorem~\ref{Kot}] 
Suppose the statement of theorem is false.
Consider a counterexample $G$ with a minimum number of edges. Let $F$ the set of all edges of  the unique perfect matching 
and $\overline{F} = E(G)\setminus F$. Consider  two cases.

$1^\circ$ {\it There exists a  bridge $a\in \overline{F}$ of the graph $G$.}

\noindent  Consider the graph $G-a$. Clearly, it has exactly two connected components. 
If any of these components has a bridge that belongs to $F$,  then the graph $G$ also has such a bridge. We obtain a contradiction.
Hence, each connected component has a second perfect matching. Clearly, then $ G $ also has second matching. This is a  contradiction.

$2^\circ$ {\it Set $\overline{F}$  contains no bridges of the graph $G$.} 

\noindent
Each vertex  is  incident to an edge of $\overline{F}$,  otherwise $G$ has a bridge that belongs to $F$. 
Furthermore, at least one vertex is incident to at least two edges of $\overline{F}$, since otherwise $G$ is an even cycle and hence, it has two perfect matchings.

Thus, $|F|<|\overline{F}|$. Since $G$ is a minimal counterexample  after deleting from $G$ any edge of $\overline{F}$ a  bridge that belongs to $F$ appears in the resulting graph.  Hence,  there are two edges $ a_1, a_2 \in \overline{F}$ such that both graphs $G-a_1$ and $G-a_2$
have the same bridge $ b \in F $.
Hence, the graph $G-b$ has  two bridges $a_1$ and $a_2$, and the graph $G-\{b,a_1,a_2\}$ has three connected components.

Returning the bridge $b$ we obtain two connected components $X, Y$ in $G-\{a_1,a_2\}$.
Since $a_1$ and $a_2$ are not bridges in $G$, we can denote them $a_1=x_1y_1$ and $a_2=x_2y_2 $, where $x_1,x_2 \in X$ $y_1,y_2  \in Y$ (see Fig.~\ref{zamena}, left).

Denote by $F_x$ the  edges of matchings $F$ lying in the component $ X $. 
Let us contract  $x_1y_1$, $Y$ and $y_2x_2$ in edge $x_1x_2$ (possibly multiple edges or a loop can appear, see Fig.~\ref{zamena}) and obtain a graph $G_x$.

{\begin{figure}[htb!]

\centering
\includegraphics[scale=0.5]{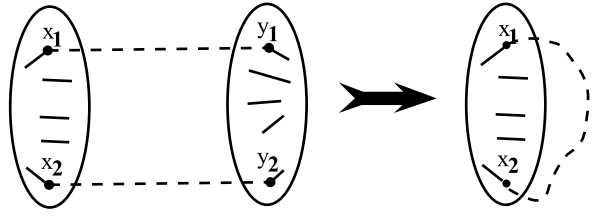}
\caption{Solid lines - edges from $F$, dotted - from $\overline{F}$.}
\label{zamena}
\end{figure}
}

Obviously, if $G_x$ has a bridge of $F_x$, then $G$ has a bridge of $F$. Hence, there is $F'_{x}$ -- a matching in $G_x$ that is different from $F_x$. If $ x_1x_2 \notin F'_{x}$, then
 there is a second matching $ F '= (F \setminus F_x) \cup F'_x$ in the graph $G$, but it's impossible. Then $F'_x$ contains an edge $x_1x_2$, 
One can similarly define the matching $F_y$ and the graph $G_y$ and prove that there exists another matching $F'_y$ in
$G_y$ that contains $y_1y_2$. Then the graph $G$ has a matching $F'= (F'_x \setminus \left\{x_1x_2 \right\})\cup (F'_y \setminus \left\{y_1, y_2 \right\})\cup \left\{x_1y_1, x_2, y_2 \right\} $, 
different from $F$, contradiction.
\end{proof}


\begin{proof}[\bf  Proof of Theorem~\ref{Ye}] 
Assume the contrary and consider a minimal counterexample $G$ (at first we minimize the number of vertices, after that the number of edges).
Clearly, $G$ is connected and  has no cut-vertex.
Vertex is called {\it monochrome}  if all edges incident to it are of the same color. And 
vertex $v$ of $G$  is called {\it cut-color} if  no connected component of $G-v$ is joined to $v$ with edges of more than one color. 
Consider several cases.

$1^\circ${\it There  is an edge $b_1 b_2$, such that the graph $G-b_1b_2$ has no monochrome vertex.}

\noindent
Since $G$ is a  minimal counterexample  there is a cut-color vertex $v$ in $G-b_1b_2$. Graph $G-b_j$ is connected, hence, $b_j \neq v$. 
Let  $X_1$ be a subgraph of the graph $G-{b_1b_2}$  induced on  the vertex $v$ together with the connected component of $ G-{b_1b_2}-v $ that contains $b_1$. Now we
add to $X_1$ the vertex $c$ and edge  $b_1c$ of the same color as $b_1b_2$ and edge  $vc$ of any unique color. Denote the obtained graph by $G_1$ (see Fig.~\ref{case1}). 
Graph $G_1$   is less than $G$, otherwise the part  $X_2$ has only one vertex and hence, the vertex $b_2$ in $G- b_1b_2$ is monochrome.

{\begin{figure}[htb!]

\centering
\includegraphics[scale=0.5]{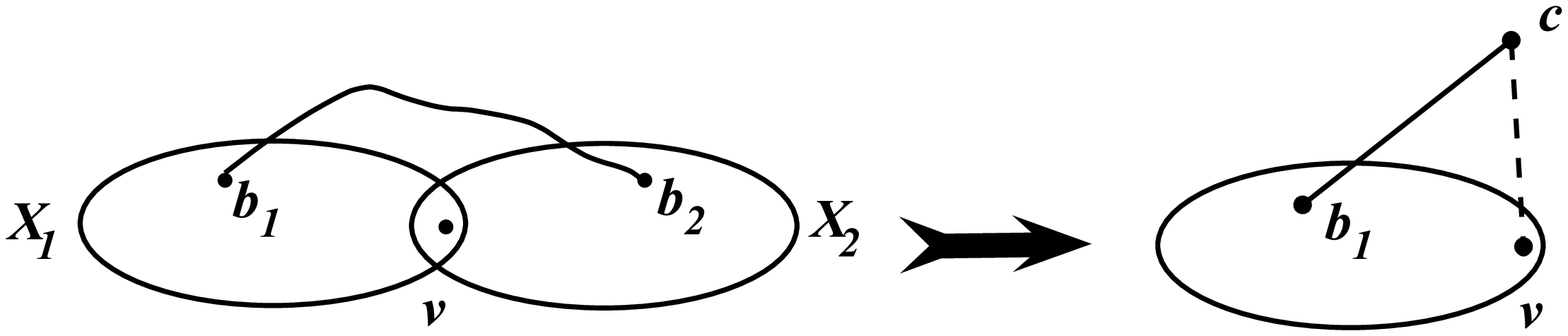}
\caption{}
\label{case1}
\end{figure}

}

The graph $G_1$ has no monochrome vertices and no cut-color vertices ($G_1$ has no cut-vertex). Hence, there is the alternating cycle in $G_1$, it must contain $C$ (otherwise, there is such a cycle in $ G $).
Consequently, there is a alternating path form $v$ to $b_1 $ in  $X_1$ with color of  last edge different from the color of edge $b_1b_2 $. Similarly, there is such a path from $v$ to $b_2 $. Two edges incident to $v$ in these paths have  different colors. Otherwise, if their colors coincide all  edges incident to  $v$ in both parts $ X_1 $ and $ X_2 $ have the same color, and therefore, $v$ is monochrome, which is impossible.  Taking these two paths and  adding to them the edge $b_1b_2$  we obtain an alternating cycle in $G$. This contradicts our assumption.

$2^\circ$ {\it There is a vertex  $c$ of degree $2$ such that the graph $G-c$ has no monochrome vertex.} 

\noindent Since $G$ is a minimal counterexample, there appears a cut-color vertex $v$ in $G-c$.
Vertex $c$ is incident to two edges  $cb_1$ and  $cb_2$ of distinct color. Let  $X_1$ be a subgraph of the graph $G-c$  induced on  the vertex $v$ together with the connected component of $ G-c-v $ that contains $b_1$.
Then we construct (similarly $1^\circ$~case) an alternating path from $v$ to $b_1$  with last edge with distinct color from $ b_1c$ in the part $X_1$ and an alternating path from $v$ to $b_2$  with last edge with distinct color from $b_2c$ in the part $X_2$.
 Glue these paths together with the edges   $b_1c$ and $cb_2$, we obtain an alternating cycle in $G$. This contradicts our assumption.

$3^\circ$ {\it There are two adjacent vertices  $c_1, c_2$ of degree 2.} 

\noindent
Let  $c_1b_1$ and $c_2b_2$ be the other edges incident to $c_1$ and $c_2$, respectively. If these two edges have different colors,
then after deleting the edge $c_1c_2$ and gluing  the vertices $c_1$ and $c_2$  we obtain a smaller counterexample. If $c_1b_1$ and $c_2b_2$ have  the same  color,  we delete the vertices  $c_1$ and $c_2$ and  add a new edge $b_1b_2$ with the same color as 
$b_1c_1$. Clearly, we obtain a smaller counterexample.

$4^\circ$ {\it Consider the remaining cases.}

\noindent Let us construct a digraph on the vertices of $G$ using his edges. We draw an arc $\vec{ab}$ if
$ab$ is an edge of $G$ and the  $b$ is a monochrome vertex of the graph $G-ab$. (Maybe the arc $\vec{ba}$ is drawn too).

Let $ x $ be the number of vertices of degree $2$ in $G$. Then the graph $G$ has at least $\frac{2x +3(v-x)}{ 2 } = 1.5v-0.5x$ edges. 
Since there is no situation  of case~$2^\circ$, at least one arc starts at each  vertex of degree $2$. Obviously, at least two arcs end at each vertex of degree 2. 
 Since  vertices of degree 2 are not adjacent, the number of arcs is  at least $ x + e \geq 1.5v +0.5 x $.

Let two arcs (corresponding to the edges $ e_1, e_2 $ of graph $ G $) end at   the vertex $ d $. Let $ E_d $ be the set of all edges incident to $ d $. Then edges of $ E_d \setminus \{e_1\} $ have the same color and edges of $ E_d \setminus \{e_2\} $ have the same color. But vertex $d$ is not monochrome in $G$, clearly, $d$ has degree 2. 
Then the sum of incoming degrees of vertices does not exceed $ 2x + (v-x) = v + x $.  But then we have $ 1.5v +0.5 x \leq v + x $, hence, $ v = x $. But we have two adjacent vertices of degree 2. We obtain  a contradiction.
\end{proof}

\end{document}